\documentclass[11pt,twoside]{article}
\usepackage{amssymb}

\evensidemargin\oddsidemargin
\newtheorem{Thm}{Theorem}

\newtheorem{Cor}[Thm]{Corollary}

\begin{document}


\thispagestyle{empty}

\begin{center}
{\Large\bf
Cohomological Splitting of Coadjoint Orbits
\medskip

}
\vskip 10mm

%
%

{\large\bf
$^1$Andr\'es Vi\~na   
}
\vskip 5mm
{\it
$^1$Departamento de F\'{\i}sica
 Universidad de Oviedo
\smallskip
 
\smallskip

email: avinae@correo.uniovi.es  
}
\end{center}

\bigskip


\begin{abstract}
\parindent0pt\noindent
 
 The rational cohomology of a coadjoint orbit ${\cal O}$ is expressed 
as tensor product of  
  the cohomology of other coadjoint orbits ${\cal O}_k$, with  
 $ \hbox{dim}\,{\cal O}_k< \hbox{dim}\,{\cal O}$.


 
\end{abstract}

%

%

\section{C-splitting of coadjoint orbits}

The purpose of this note is to express the rational cohomology of a given 
coadjoint orbit of a compact Lie group in terms of the cohomology of ``smaller"
coadjoint orbits. Our result is based upon two facts: The coadjoint orbit hierarchy,
and the cohomological splitting of certain Hamiltonian bundles.

\smallskip

{\it The coadjoint orbit hierarchy.}

Let $G$ be a compact and connected Lie group. We consider the coadjoint action of $G$
on ${\frak g}^*\,$. By
$X_A$ is denoted the vector field on ${\frak g}^*$ generated by $A\in{\frak g}$. If
 $\mu\in{\frak g}^*$, we denote
by ${\cal O}=G\cdot \mu$  the coadjoint orbit of $\mu$. Then
${\cal O}=G/G_{\mu}$, where $G_{\mu}$ is the subgroup of isotropy of $\mu$.
The manifold ${\cal O}$ possesses a natural symplectic structure defined by the
$2$-form $\omega$, with $\omega_{\nu}(X_A,X_B)=\nu([A,B])$, for 
any $\nu\in{\cal O}$ \cite{aK76}.
If $l_g$ denotes the left multiplication by $g\in G$; that is,
$l_g:\nu\in{\cal O}\mapsto g\cdot\nu\in{\cal O}$, then $l_g^*\omega=\omega$. Moreover
$\iota_{X_A}\omega=dh_A$, with $h_A$ the function on ${\cal O}$ defined by
$h_A(\nu)=\nu(A)$. Therefore the action of $G$ on ${\cal O}$ is Hamiltonian; that is,
$G$ is a subgroup of the group $\hbox{Ham}({\cal O})$ of Hamiltonian
 symplectomorphisms \cite{lP01}
of ${\cal O}$. And using Morse theory one can prove that  
  ${\cal O}$ is simply-connected \cite{GLS}.

On the other hand, if $\mu_1,\mu_2\in{\frak g}^*$ and $G_1:=G_{\mu_1}\subset G_{\mu_2}=:G_2$,
then the  orbits ${\cal O}_j=G\cdot\mu_j$,$\;j=1,2$ are in the following hierarchy:
There is a symplectic fibration of ${\cal O}_1$ over ${\cal O}_2$. In fact
${\cal O}_1=G\times_{G_2}(G_2\cdot \mu_1)$. So ${\cal O}_1$ is a fiber bundle over 
$G/G_2$ with fiber the orbit of $G_2\cdot\mu_1$ of $G_2$. Thus the fiber is in turn
 a symplectic manifold, and on it   the group $G_2$ acts as a group of Hamiltonian
 symplectomorphisms, if $G_2$ is connected (for details see \cite{GLS}).

  \smallskip

{\it Cohomological splitting of Hamiltonian bundles.}

Let $P\rightarrow B$ be a fiber bundle, with fiber a symplectic manifold $M$. This
bundle is said to be Hamiltonian if its structural group reduces to the group 
$\hbox{Ham}(M)$ of Hamiltonian symplectomorphisms of $M$ \cite{Mc-S}.
Lemma 4.11 of \cite{L-M} states that the rational cohomology of
any Hamiltonian fiber bundle $M\rightarrow P\rightarrow {\cal O}$, 
whose base is a coadjoint orbit, splits additively as the tensor product
of the cohomology of the fiber by the one of ${\cal O}$; that is, 
$H^*(P)\simeq H^*(M )  \otimes H^*({\cal O})$.

\smallskip

If we apply the result of Lalonde and McDuff to our Hamiltonian fibration
 $$G_2\cdot\mu_1\rightarrow {\cal O}_1\rightarrow {\cal O}_2,$$
we obtain an additive isomorphism
$$H^*({\cal O}_1,\,{\Bbb Q})\simeq H^*({\cal O}_2,\,{\Bbb Q})\otimes H^*(G_2/G_1,\,{\Bbb Q}),$$
in other words
\begin{Thm}\label{teoremp}
If $G_1\subset G_2$ are stabilizers of the coadjoint action of the compact,
 connected Lie group $G$ and $G_2$ is connected, then
there is an additive isomorphism
 $$H^*(G/G_1 ,\,{\Bbb Q} )\simeq H^*(G/G_2,\,{\Bbb Q})\otimes H^*(G_2/G_1,\,{\Bbb Q}).$$ 
\end{Thm}  

\begin{Cor}\label{Torre}
If $\mu_1,\dots,\mu_k$ are points of ${\frak g}^*$, such
that
$$ G_{\mu_1} \subset G_{\mu_2}\subset\dots \subset G_{\mu_k}\ne G,$$
 and the $G_{\mu_j}$ are connected, then
$$H^*(G/G_{\mu_1})\simeq H^*(G/G_{\mu_k})\otimes\bigotimes_{j=2}^k H^*(G_{\mu_j}/G_{\mu_{j-1}}).$$
\end{Cor}
This formula expresses the rational cohomology of the orbit ${\cal O}_1=G/G_{\mu_1}$ in terms of
the cohomology of  orbits whose dimensions are  less than $\hbox{dim}\,{\cal O}_1$.


\section{Cohomological splitting of flag manifolds}

A partition ${\frak p}$ of an integer $n$ is an unordered sequence $i_1,\dots,i_s$ of positive
integers with sum $n$. This partition of $n$ determines  the 
subgroup 
$$G_{\frak p}:=U(i_1)\times\dots\times U(i_s)$$
 of $U(n).$
  Moreover this subgroup
is  a stabilizer for the coadjoint action of $U(n)$.
The partitions $(11\dots 1)$, $(1\dots 12),$ $\dots,(1n-1)$ of $n$ determine a tower
of subgroups 
$$G_1\subset G_2\subset\dots\subset G_{n-1}$$
of $U(n)$. The quotient $U(n)/G_1$ is the flag manifold ${\cal F}_n$, i.e. the manifold
of complete flags  in ${\Bbb C}^n$, and
$G_j/G_{j-1}\simeq U(j)/(U(1)\times U(j-1))={\Bbb C}P^{j-1}$. From Corollary \ref{Torre}
we deduce
\begin{Cor}\label{Flag}
If ${\cal F}_n$ denotes the flag manifold in ${\Bbb C}^n$, then
$$H^*({\cal F}_n,{\Bbb Q})\simeq\bigotimes_{j=1}^{n-1}H^*({\Bbb C}P^{j},{\Bbb Q}).$$
\end{Cor}

As particular case we consider the group $G:=U(4)$ and its subgroups
$$G_1:=U(1)\times\dots\times U(1)\subset U(2)\times U(2)=:G_2.$$
Then by Theorem \ref{teoremp} 
\begin{equation}\label{flag4}
H^*({\cal F}_4)\simeq H^*(G_{2,2}({\Bbb C})) \otimes H^*(G_2/G_1)
\simeq H^*(G_{2,2}({\Bbb C})) \otimes H^*({\Bbb C}P^1)\otimes H^*({\Bbb C}P^1),
\end{equation}
where $G_{2,2}({\Bbb C})$ is the corresponding Grassman manifold in ${\Bbb C}^4$. 
So by Corollary \ref{Flag}
\begin{equation}\label{flag4I}
H^*(G_{2,2}({\Bbb C}))\otimes H^*({\Bbb C}P^1)\otimes H^*({\Bbb C}P^1)
\simeq H^*({\Bbb C}P^1)\otimes H^*({\Bbb C}P^2)\otimes H^*({\Bbb C}P^3).
\end{equation}
 The existence of this isomorphism can be checked directly.
The cohomology $H^*(G_{2,2})$ is generated by $\{c_1,c_2\}$,
where $c_i$ is the corresponding Chern class of the $2$-plane universal bundle over
$G_{2,2}({\Bbb C})$ (see \cite{Dold}). Moreover $c_1,c_2$ are algebraically independent
up to dimension $4$. So $\hbox{dim}\, H^{4}(G_{2,2})=2$ and 
$\hbox{dim}\, H^{2j}(G_{2,2})=1$, for $j\ne 2,\;0\leq j\leq 4$. 
Therefore it is possible to identifify
the graded vector spaces $H^*(G_{2,2})\otimes H^*({\Bbb C}P^1)$ 
and $H^*({\Bbb C}P^3)\otimes H^*({\Bbb C}P^2)$. This identification allows
 us to construct the isomorphism (\ref{flag4I}).

In general, a partition ${\frak p}$ of $n$ determines  
 the manifold of partial flags 
${\cal F}_{\frak p}=U(n)/G_{\frak p}$. 
The following corollary is a consequence of
  Theorem \ref{teoremp}
 
\begin{Cor}
If ${\frak p}=\{i_1,\dots,i_s\}$ and ${\frak p'}=\{j_1,\dots,j_r\}$ are partitions of $n$ with
$r<s$ and  $G_{\frak p}\subset  G_{\frak p'} $, then
$$ H^*({\cal F}_{\frak p},{\Bbb Q})= H^*({\cal F}_{\frak p'},{\Bbb Q})
\otimes H^*(G_{\frak p'}/G_{\frak p},{\Bbb Q}).$$
\end{Cor}

\section*{Acknowledgments}

The author was partially supported by Universidad de Oviedo, grant NP-01-514-4




\begin{thebibliography}{99}   
 
    


 \bibitem{Dold}
 { A. Dold}, 
 {\em Lectures on Algebraic Topology}, 
  Springer, 
   Berlin. 
  (1980) 




\bibitem{GLS}
{ V. Guillemin,  E. Lerman, S. Sternberg},  
 {\em Symplectic Fibrations and Multiplicity Diagrams},
Cambridge U. P.,
Cambridge. (1996)










\bibitem{aK76}
{ A.\,A. Kirilov}
{\em  Elements of the Theory of Representations}, 
  Springer-Verlag, 
  Berlin. 
(1976)



    


    
 \bibitem{L-M}
F. Lalonde, D. McDuff, {\em Symplectic Structures on Fiber Bundles}
 Topology 42(2) (2002) 309-347 
 



      
 
 
       
 
 





\bibitem{Mc-S}
  {D. McDuff, D. Salamon},
{\em Introduction to Symplectic Topology},
Clarenton Press,
 Oxford. 
(1998)
 


\bibitem{lP01}
{ L. Polterovich},
{\em The Geometry of the Group of Symplectic Diffeomorphisms},
Birkh\"auser,
Basel. (2001)



  

   
  
   

 

 
\end{thebibliography}
\end{document}